\documentclass[11pt]{article}
\usepackage{amssymb,amsfonts,amsmath,latexsym,epsf,tikz,url}
\newtheorem{theorem}{Theorem}[section]
\newtheorem{proposition}[theorem]{Proposition}

\newcommand{\proof}{\noindent{\bf Proof.\ }}
\newcommand{\qed}{\hfill $\square$\medskip}

\begin{document}

\title{ Total dominator chromatic number of $k$-subdivision of graphs}

\author{Saeid Alikhani$^{}$\footnote{Corresponding author} \and Nima Ghanbari \and Samaneh Soltani  }

\date{\today}

\maketitle

\begin{center}
   Department of Mathematics, Yazd University, 89195-741, Yazd, Iran\\
\bigskip{\tt alikhani@yazd.ac.ir\\  n.ghanbari.math@gmail.com\\ s.soltani1979@gmil.com}\\
\end{center}


\begin{abstract}

Let $G$ be a simple graph. A total dominator coloring of $G$, is a proper coloring of the vertices
of $G$ in which each vertex of the graph is adjacent to every vertex of some color class.
The total dominator chromatic (TDC) number $\chi_d^t(G)$ of $G$, is the minimum number of colors
among all total dominator coloring of $G$. For any $k \in \mathbb{N}$, the $k$-subdivision of $G$ is a simple graph $G^{\frac{1}{k}}$ which is constructed by replacing each edge of $G$ with a path of length $k$. 
In this paper, we study the total dominator chromatic number of $k$-subdivision of $G$. 

\end{abstract}


\medskip
\noindent{\bf AMS Subj.\ Class.:} 05C15, 05C69

\section{Introduction}

In this paper, we are concerned with simple finite graphs. Let $G=(V,E)$ be such a graph and $\lambda \in \mathbb{N}$. A mapping $f : V (G)\longrightarrow \{1, 2,...,\lambda\}$ is
called a \textit{$\lambda$-proper  coloring} of $G$, if $f(u) \neq f(v)$ whenever the vertices $u$ and $v$ are adjacent
in $G$. A \textit{color class}  of this coloring, is a set consisting of all those vertices
assigned the same color. If $f$ is a proper coloring of $G$ with the coloring classes $V_1, V_2,..., V_{\lambda}$ such
that every vertex in $V_i$ has color $i$, then sometimes write simply $f = (V_1,V_2,...,V_{\lambda})$.  The \textit{chromatic number } $\chi(G)$ of $G$ is
the minimum number of colors needed in a proper coloring of a graph.
The concept of a graph coloring and chromatic number is very well-studied in graph theory.
The total dominator coloring, abbreviated TD-coloring studied in  \cite{Adel,Adel2,Vij1,Vij2}. Let $G$ be a graph with no
isolated vertex, the \textit{ total dominator coloring}  is a proper coloring of $G$ in which each vertex of the graph is adjacent
to every vertex of some (other) color class. The \textit{total dominator chromatic number}, abbreviated TDC-number, $\chi_d^t(G)$ of $G$ is the minimum number of color classes in a TD-coloring of $G$.  
A set $D\subset V$ is a total dominating set if every vertex of $V$ is adjacent to some vertices of $D$. The total dominating number $\gamma_t(G)$ is the minimum cardinality of a total dominating set in $G$.

 Computation of the TDC-number is NP-complete (\cite{Adel}).
The TDC-number of some graphs, such as paths, cycles, wheels and the complement of paths and cycles has computed in
\cite{Adel}. Also Henning in \cite{GCOM} established the  lower and the upper bounds on the TDC-number
of a graph in terms of its total domination number $\gamma_t(G)$. He has shown that,  every
graph $G$ with no isolated vertex satisfies $\gamma_t(G) \leq \chi_d^t (G)\leq \gamma_t(G) + \chi(G)$.
The properties of TD-colorings in trees has studied in \cite{GCOM,Adel}. Trees $T$ with $\gamma_t(T) =\chi_d^t(T)$ has characterized
in \cite{GCOM}. We have examined the effects on $\chi_d^t(G)$ when $G$ is modified by operations on vertex and edge of $G$ in \cite{nima2}.

The \textit{ $k$-subdivision} of $G$, denoted by $G^{\frac{1}{k}}$, is constructed by replacing each edge $v_iv_j$ of $G$ with a path of length $k$., say $P_{\{v_i,v_j\}}$. These $k$-paths are called \textit{superedges}, any new vertex is an internal vertex, and is denoted by $x^{\{v_i,v_j\}}_l$ if it belongs to the superedge $P_{\{v_i,v_j\}}$, $i<j$ with  distance $l$ from the vertex $v_i$, where $l \in \{1, 2, \ldots , k-1\}$.  Note that for $k = 1$, we have $G^{1/1}= G^1 = G$, and if the graph $G$ has $v$ vertices and $e$ edges, then the graph $G^{\frac{1}{k}}$ has $v+(k-1)e$ vertices and $ke$ edges.

\medskip

In this paper, we study the TDC-number of the  $k$-subdivision of  a graph. 

\section{Main results}

We start by proposing  upper and lower  bounds for the TDC-number of the  $k$-subdivision of  a graph. 
First we need the TDC-number  of path graph. Note that the value of  TDC-number  of paths and cycles which have computed in \cite{Adel} are lower and upper bounds for
 $\chi_d^t(P_n)$,  $\chi_d^t(C_n)$  and are  not
the exact value. For example by formula in \cite{Adel}, 
 	$\chi_ d^t (P_{60})=40$ which is not true and the correct value is $32$ which can obtain by the following theorem.

  \begin{theorem}\label{newpath}
  	If  $P_n$ is  the path graph of order $n\geq 8$, then
 \[
 	\chi_d^t(P_n)=\left\{
  	\begin{array}{ll}
  	{\displaystyle
  		2k+2}&
  		\quad\mbox{if $n=4k$, }\\[15pt]
  		{\displaystyle
  			2k+3}&
  			\quad\mbox{if $n=4k+1$,}\\[15pt]
  			{\displaystyle
  				2k+4}&
  				\quad\mbox{if $n=4k+2$, $n=4k+3$.}\\[15pt]
  				  					\end{array}
  					\right.	
  					\]
  Also $ \chi _ d^t (P_3)=2$, $ \chi _ d^t (P_4)= 3$, $\chi _ d^t (P_5)= \chi _ d^t (P_6)=4$ and $ \chi _ d^t (P_7)=5$.
  					
  				\end{theorem}
  				
\begin{figure}
  					\begin{center}
  					\includegraphics[width=1.7in]{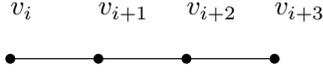}
  						\caption{Four consecutive vertices of the Path graph $P_n$.}
  					\label{Four}
  				\end{center}
  				\end{figure}

  		\proof 	It is easy to show that $ \chi _ d^t (P_3)=2$, $ \chi _ d^t (P_4)= \chi _ d^t (P_5)=3$, $ \chi _ d^t (P_6)=4$ and $ \chi _ d^t (P_7)=5$. Now let $n\geq 8$. First we show that for each four consecutive vertices we have to use at least two new colors. Consider Figure \ref{Four}. We have two cases. If we give an old color to $v_{i+1}$, then we need to give a new color to $v_{i+2}$ and $v_{i+3}$ to have a TD-coloring. Also if we give a new color to $v_{i+1}$, then we have to give a new color to $v_{i+2}$ or $v_{i}$ to have a TD-coloring. So we need at least two new colors in every four consecutive vertices. 
  				
  				Suppose that $n=4k$, for some $k\in \mathbb{N}$. We give a TD-coloring for  the path $P_{4k}$ which use only two new colors in every four consecutive vertices. Define a function $f_0$ on the vertices of $P_{4k}$, i.e., $V(P_{4k})$ such that for any vertex  $v_i$, 
  	 \begin{center}
  					 $f_0(v_i)=\left\{
  					\begin{array}{ll}
  					{\displaystyle
  						1} &
  						\quad\mbox{if $i=1+4k$,}\\[15pt]
  					{\displaystyle
  						2} &
  						\quad\mbox{if $i=4k$,}
  					\end{array}
  					\right.$
  \end{center}
  	 				and for any $v_i$ , $i \neq 4k$ and $i\neq 4k+1$, $f_0(v_i)$ is a new number. Then $f_0$ is a TD-coloring  of $P_{4k}$ with the minimum number $2k+2$.
  				
  		If $n=4k+1$, for some $k\in \mathbb{N}$, then first color the $4k-4$ vertices  using  $f_0$. Now for the rest of vertices define $f_1(v_{4k-3})=1$, $f_1(v_{4k-2})=2k+1$, $f_1(v_{4k-1})=2k+2$, $f_1(v_{4k})=2k+3$ and $f_1(v_{4k+1})=2$. Since for every five consecutive vertices we have to use at least three new colors, so  $f_1$ is a TD-coloring  of $P_{4k+1}$ with the minimum number $2k+3$.

  				If  $n=4k+2$, for some $k\in \mathbb{N}$, then first color  the $4k-4$ vertices using  $f_0$. Now for the rest of vertices define $f_2(v_{4k-3})=1$, $f_2(v_{4k-2})=2k+1$, $f_2(v_{4k-1})=2k+2$, $f_2(v_{4k})=2k+3$, $f_2(v_{4k+1})=2k+4$ and $f_2(v_{4k+2})=2$. Since  for every six consecutive vertices we have to use at least four  new colors, so  $f_2$ is a TD-coloring  of $P_{4k+2}$ with the minimum number $2k+4$.

  				If  $n=4k+3$, for some $k\in \mathbb{N}$, then first color  the $4k-4$ vertices using $f_0$. Now for the rest of vertices define $f_3(v_{4k-3})=1$, $f_3(v_{4k-2})=2k+1$, $f_3(v_{4k-1})=2k+2$, $f_3(v_{4k})=2$, $f_3(v_{4k+1})=2k+3$, $f_3(v_{4k+2})=2k+4$ and $f_3(v_{4k+2})=2$. Then $f_3$ is a TD-coloring of $P_{4k+2}$ with the minimum number $2k+4$. Therefore we have the result. \qed

 \begin{theorem}\label{ubboun2}
 If  $G$ is a connected graph with $m$ edges and $k\geq 2$, then 
 \begin{equation*}
 \chi_d^t(P_{k+1}) \leq \chi_d^t(G^{\frac{1}{k}})\leq (m-1) \chi_d^t(P_{k})+\chi_d^t(P_{k+1}).
 \end{equation*}
 \end{theorem}
 \proof  For the right inequality, let $e=uu_1$ be an arbitrary edge of $G$. This edge is replaced with   the superedge $P^{\{u,u_1\}}$  in $G^{\frac{1}{k}}$, with vertices $\{u, x^{\{u,u_1\}}_1, \ldots , x^{\{u,u_1\}}_{k-1}, u_1\}$. We color this superedge with  $\chi_d^t(P_{k+1})$ colors as a total dominator coloring of $P_{k+1}$ (Theorem \ref{newpath}). If $N_G(u)=\{u_1, \ldots , u_s\}$ then we color  the vertices of paths $P^{\{u,u_i\}}$ such that
 \begin{enumerate}
 \item The color of $u$ in $P^{\{u,u_i\}}$, for any $2\leq i \leq s$, is the same as color $u$ in total dominator coloring of $P^{\{u,u_1\}}$.  
 \item The superedges $P^{\{u,u_i\}}$, for any $1\leq i \leq s$, have been colored as a total dominator coloring of $P_{k+1}$ such that for any $i,i'\in \{1, \ldots , s\}$, (note that $c(y)$ is the color of vertex $y$ in our coloring)
 \begin{equation*}
 \left(\bigcup_{j=1}^{k-1}c(x^{\{u,u_i\}}_{j}) \cup c(u_i) \right)\bigcap  \left(\bigcup_{j=1}^{k-1}c(x^{\{u,u_{i'}\}}_{j}) \cup c(u_{i'}) \right)=\emptyset,~\text{where} ~i\neq i'.
 \end{equation*}
 \end{enumerate}
  Thus we need   at most $ (s-1) \chi_d^t(P_{k})+\chi_d^t(P_{k+1})$ colors for such coloring of vertices of superedges $P^{\{u,u_i\}}$,  $1\leq i \leq s$. Note that we need at most $\chi_d^t(P_{k})$ new colors for a TD-coloring of $P^{\{u,u_i\}}$,  since the vertex $u$ has been colored in all  superedges $P^{\{u,u_i\}}$,  $1\leq i \leq s$. We do not use the colors used for superedges $P^{\{u,u_i\}}$,  $1\leq i \leq s$, any more. In the next step, we consider that superedges in $G^{\frac{1}{k}}$ which are replaced instead of incident edges to $u_i$'s in $G$, and have not been colored in the prior step. Now we color the vertices of these superedges as a TD-coloring of $P_{k}$, such that the  vertices $u_2, \ldots , u_s$ have been colored in prior step, and the pairwise intersection of the set of colors used for coloring of vertices of these superedges is the empty set. We continue this process to color all vertices  of $G^{\frac{1}{k}}$.  This coloring is   a TD-coloring of $G^{\frac{1}{k}}$, because every superedge have been colored with distinct color set, except the end vertices of the superedges, possibly. Finally, since we used at most  $(m-1) \chi_d^t(P_{k})+\chi_d^t(P_{k+1})$ colors, the right inequality follows.
 
  For the left inequality, if $G$ is a path then the result is true. So we suppose that $G$ is a connected graph which is not a path. Let  $P^{\{v,w\}}$ be an arbitrary superedge of $G^{\frac{1}{k}}$ with vertex set  $\{v, x^{\{v,w\}}_1, \ldots , x^{\{v,w\}}_{k-1}, w\}$. Since $G$ is not a path, so at least one of $v$ and $w$ is adjacent to some vertices of $G^{\frac{1}{k}}$ except $x^{\{v,w\}}_1$ and $x^{\{v,w\}}_{k-1}$, respectively. Let $c$ be a total dominator coloring of  $G^{\frac{1}{k}}$. The two following cases can be occured: either the restriction of $c$ to vertices of $P^{\{v,w\}}$ is a total dominator coloring and so we have the result, or not. If the restriction of $c$ to vertices of $P^{\{v,w\}}$ is not a total dominator coloring then since $c$ is a total dominator coloring  of $G^{\frac{1}{k}}$ we conclude that at least one of vertices $v$ and $w$, as the vertices of the induced subgraph $P^{\{v,w\}}$, are not adjacent to every vertex of some color class. Without loss of generality we assume that the vertex $v$, as the vertex of the induced subgraph $P^{\{v,w\}}$,  is not  adjacent to every vertex of some color class. But $c$ is a total dominator coloring of $G^{\frac{1}{k}}$ so the vertex $v$ is adjacent to every vertex of some color class, as the vertex of  $G^{\frac{1}{k}}$. Hence there is a new color for an adjacent vertex of $v$, except the vertex $x^{\{v,w\}}_1$. Thus if we use this new color for the vertex  $x^{\{v,w\}}_1$ and consider the restriction of $c$ for the remaining vertices of superedge $P^{\{v,w\}}$, then $P^{\{v,w\}}$ has a total dominator coloring. Therefore the total coloring $c$ has at least  $\chi_d^t(P_{k+1}) $ colors.\qed
 
   By the following Proposition we  show that the upper bound of  $\chi_d^t(G^{\frac{1}{k}})$ in Theorem \ref{ubboun2} is sharp for $G=K_{1,n}$ and $k=3$.
  \begin{proposition}
 For every $n\geq 3$,  $\chi_d^t(K_{1,n}^{\frac{1}{3}}) = 2n+1$.
  \end{proposition}
\proof  Let $p_1, \ldots , p_n$ be the pendant vertices of $K_{1,n}^{\frac{1}{3}}$. The adjacent vertex to $p_i$ is denoted by $q_i$, and the adjacent vertex to $q_i$ of degree 2 is denoted by $w_i$ for any $1\leq i \leq n$. The center of $K_{1,n}^{\frac{1}{4}}$ is denoted by $v$.  Since  the vertex  $q_i$ is the only vertex adjacent to $p_i$, so the color of $q_i$ should not be used for any other vertices of graph, where $1\leq i \leq n$. Thus we color the vertices $q_1, \ldots , q_n$ with colors $1, \ldots , n$, respectively, and do not use these colors any more. For every $1\leq i \leq n$, the vertex $q_i$ is adjacent to $p_i$ and $w_i$, thus we need a new color for at least one of $p_i$ and $w_i$. Here we consider the three following cases:
  \begin{itemize}
  \item[Case 1)] The vertices $p_1, \ldots , p_n$ have been colored  with 
  colors $n+1, n+2, \ldots , 2n$, respectively, and these colors do not use any more. In this case we must color the vertices $w_1, \ldots , w_n$ with a new color, say color $2n+1$. Since we need a proper coloring, we color the vertex $v$ with color $2n+2$. So it can be seen that we have a TD-coloring of $K_{1,n}^{\frac{1}{3}}$ with $2n+2$ colors.
  \item[Case 2)] The vertices $w_1, \ldots , w_n$ have been colored  with colors $n+1, n+2, \ldots , 2n$, respectively, and these colors do not use any more. Since we need a proper coloring, we must color the vertex $v$ with a new color. If we color all vertices $v, p_1, \ldots , p_n$ with color $2n+1$ then we have a TD-coloring of $K_{1,n}^{\frac{1}{3}}$ with $2n+1$ colors.
  \item[Case 3)] The vertices $w_i,p_i$ have been colored  with color $n+i$, for any $i$ $1\leq i \leq n$ and these colors do not use any more. Thus we need a new color for the vertex $v$, say  color $2n+1$. Now we have a TD-coloring of $K_{1,n}^{\frac{1}{3}}$ with $2n+1$ colors.
  \end{itemize}

Therefore   $\chi_d^t(K_{1,n}^{\frac{1}{3}}) = 2n+1$, by Cases 1, 2 and 3.  \qed

  Here we improve the  lower bound of Theorem \ref{ubboun2} for  $k\geq 9$.
  \begin{theorem}\label{lbound}
 If  $G$ is a connected graph with $m$ edges and $k\geq 9$, then 
 \begin{equation*}
 m(\chi_d^t(P_{k-1})-2) +2 \leq \chi_d^t(G^{\frac{1}{k}}).
 \end{equation*}
 \end{theorem}
\proof Let $e=vw$ be an edge of $G$. We consider the superedge $P^{\{v,w\}}$ with vertex set $\{v, x^{\{v,w\}}_1, \ldots , x^{\{v,w\}}_{k-1}, w\}$. It is clear that $P^{\{v,w\}}\setminus \{v,w\}$ is the path graph $P_{k-1}$. Since we use  repetitious colors for the vertices $x^{\{v,w\}}_1$ and $x^{\{v,w\}}_{k-1}$  in the  TD-coloring of paths, so we need at least $\chi_d^t(P_{k-1})-2$ colors for each superedges, and hence  the result follows.\qed

   \begin{theorem}\label{lower1}
 If  $G$ is  a connected graph with $m$ edges and $k\geq 9$, then  
 \begin{equation*}
 \chi_d^t(G^{\frac{1}{k}})\geq \left\{
 \begin{array}{ll}
 \frac{mk}{2}+2 & k\equiv 0 ~(\text{mod} ~4)\\
 m(\frac{k-1}{2})+2 &k\equiv 1 ~(\text{mod} ~4)\\
 m(\frac{k-2}{2}+1)+2 & k\equiv 2 ~(\text{mod} ~4)\\
 m(\frac{k-3}{2}+2)+2 & k\equiv 3 ~(\text{mod} ~4).
 \end{array}\right.
 \end{equation*}
  \end{theorem}
   \proof 
 It follows by Theorems \ref{newpath} and \ref{lbound} .\qed

\begin{theorem}\label{ubound}
 If  $G$ is a connected graph with $m$ edges and $k\geq 7$, then 
 \begin{equation*}
  \chi_d^t(G^{\frac{1}{k}})\leq m(\chi_d^t(P_{k+1})-2) +2.
 \end{equation*}
 \end{theorem}
\proof As we have seen  in TD-coloring  of paths, we can use the same color for the pendant vertices. So we give the  color $1$ or $2$ to all the vertices belong to $G$ and we color other vertices of any superedges with $\chi_d^t(P_{k-1})-2$ colors. This is a TD-coloring for $G^{\frac{1}{k}}$ and hence  the result follows.\qed

\begin{theorem}
If  $G$ is  a connected graph with $m$ edges and $k\geq 9$, then  
 \begin{equation*}
 \chi_d^t(G^{\frac{1}{k}})\leq \left\{
 \begin{array}{ll}
 m(\frac{k}{2}+1)+2 & k\equiv 0 ~(\text{mod} ~4)\\
 m(\frac{k-1}{2}+2)+2 &k\equiv 1 ~(\text{mod} ~4)\\
 m(\frac{k-2}{2}+2)+2 & k\equiv 2 ~(\text{mod} ~4)\\
 m(\frac{k-3}{2}+2)+2 & k\equiv 3 ~(\text{mod} ~4).
 \end{array}\right.
 \end{equation*}
  \end{theorem}
   \proof 
 It follows by Theorems \ref{newpath} and \ref{ubound} .\qed

Iradmusa \cite{Moharram}, showed that if $G$ is a connected graph and $k$ is a positive integer greater than one, then at most three colors are enough to achieve
a proper coloring of $G^{\frac{1}{k}}$.
 
 \begin{theorem}
 If  $G$ is  a connected graph and $k$ be a positive integer number greater than one, then $\gamma_t(G^{\frac{1}{k}}) \leq \chi_d^t (G^{\frac{1}{k}})\leq \gamma_t(G^{\frac{1}{k}}) + 2$.
 \end{theorem}
\proof If $G$ is a connected graph and $k$ is a positive integer greater than one, then at most three colors are enough to achieve
a proper coloring of $G^{\frac{1}{k}}$. In continue we want to show that there is no graph with $\chi_d^t (G^{\frac{1}{k}})= \gamma_t(G^{\frac{1}{k}}) + 3$. 
 Let $\gamma_t(G^{\frac{1}{k}})=s$ and $\Gamma =\{y_1, \ldots , y_s\}$ be the total dominating set of $G^{\frac{1}{k}}$. 
We color the vertices $y_1, \ldots , y_s$ with colors $1, \ldots , s$, respectively,  and do not use these colors any more. Let $x$ be an arbitrary uncolored vertices of $G^{\frac{1}{k}}$ which is adjacent to one of  $y_1, \ldots , y_s$. 
Since $k\geq 2$, so the set of all uncolored adjacent vertices to $x$, denote by $N_{un-col}(x)$, make an independent set of vertices, since otherwise we find  a triangle, i.e., $K_3$, in $G^{\frac{1}{k}}$, which is a contradiction (note that if $G$ is a graph with girth $g$, then $G^{\frac{1}{k}}$ is  graph with girth $gk$). If $N_{un-col}(x)=\emptyset$, then we color $x$ with color $s+1$. If $N_{un-col}(x)\neq \emptyset$, then we suppose that $N_{un-col}(x)=\{x_1, \ldots , x_t \}$ for some $t\geq 1$. 
 Since $k\geq 2$ and $\{x_1, \ldots , x_t\} \cap \{y_1, \ldots , y_s\}=\emptyset$,  so each of vertices $x_i$ has no uncolored adjacent vertex in  $G^{\frac{1}{k}}$, except $x$,  because  otherwise $\{y_1, \ldots , y_s\}$ is not a total dominating set. In this case we color the vertex $x$ with color $s+1$, and all vertices $x_1, \ldots , x_t$ with color $s+2$. It can be seen that we have a TD-coloring of   $G^{\frac{1}{k}}$.  By the above argument, we conclude that  $\chi_d^t (G^{\frac{1}{k}})\leq \gamma_t(G^{\frac{1}{k}}) + 2$, for any $k\geq 2$. \qed

 
 The following proposition gives the exact value TDC-number  of $4$-subdivision of 
 stars graph.  
  \begin{proposition}
  For every $n\geq 3$, $\chi_d^t(K_{1,n}^{\frac{1}{4}}) = 2n+2$.
  \end{proposition}
  \proof  Let $p_1, \ldots , p_n$ be the pendant vertices of $K_{1,n}^{\frac{1}{4}}$. The adjacent vertex to $p_i$ is denoted by $q_i$, and the adjacent vertex to $q_i$ of degree 2 is denoted by $w_i$ for any $1\leq i \leq n$. If the center of $K_{1,n}^{\frac{1}{4}}$ is denoted by $v$ then the adjacent vertices of $v$ are denoted by $z_1, \ldots , z_n$ where $z_i$ and $w_i$ are adjacent for any $1\leq i \leq n$.  Since $p_i$ is a vertex of degree one, for any $1\leq i \leq n$, and $p_i$ must be adjacent to all vertices of a color class, so we must color the vertices $q_1, \ldots , q_n$ with different colors, say $1, \ldots , n$, respectively, and we should not use these colors any more. Now we consider the vertices  $q_1, \ldots , q_n$. These vertices have exactly two adjacent vertices, $w_i$'s and $p_i$'s. To have a TD-coloring, at least one of these two vertices,   $w_i$ and $p_i$, must be colored with a new color, for every $1\leq i \leq n$, such that these new colors do not use any more.  Hence the three following cases can be occured:
  \begin{itemize}
  \item[Case 1)] The vertices $w_1, \ldots , w_n$ have been colored  with colors $n+1, n+2, \ldots , 2n$, respectively, and these colors do not use any more. It is clear that we need at least two colors for coloring of vertices $v, z_1, \ldots , z_n$, since the vertex $v$ is adjacent to vertices $z_1, \ldots , z_n$. If we label all vertices $z_1, \ldots , z_n$ with color $2n+1$, and the vertices $v, p_1, \ldots , p_n$ with color $2n+2$, then it can be seen that we have a TD-coloring of $K_{1,n}^{\frac{1}{4}}$ with $2n+2$ colors.
  \item[Case 2)] The vertices $p_1, \ldots , p_n$ have been colored  with colors $n+1, n+2, \ldots , 2n$, respectively, and these colors do not use any more. In this case we consider the vertices $z_1, \ldots , z_n$. Each of $z_i$, $1\leq i \leq n$, has two adjacent vertices, $v$ and $w_i$. Since $z_i$ must be adjacent to all vertices of some color class, so at least one of $w_i$ and $v$ must be colored with a new color such that this color do not use any more. Since $v$ is the common neighbors of all  vertices $z_1, \ldots , z_n$, and we want to use the minimum number of colors, so we color  the vertex $v$ with color $2n+1$, and do not use this color any more. Since our coloring must be proper, we need at least two new colors for coloring of the remaining uncolored vertices: one color, say $2n+2$, for vertices $z_1, \ldots , z_n$, and another color for coloring of all vertices $w_1, \ldots , w_n$, say color $2n+3$. Now we have a TD-coloring of $K_{1,n}^{\frac{1}{4}}$ with $2n+3$ colors.
  \item[Case 3)] The vertices $w_i,p_i$ have been colored  with color $n+i$, for any $i$ ($1\leq i \leq n$) and these colors do not use any more. Since our coloring is proper so we need at least two new colors  for coloring of vertices $v$ and $z_1, \ldots , z_n$. If we color the vertex $v$ with color $2n+1$, and all vertices $z_1, \ldots , z_n$ with color $2n+2$, then it can be seen that we have a TD-coloring of $K_{1,n}^{\frac{1}{4}}$ with $2n+2$ colors.
  \end{itemize}
Therefore   $\chi_d^t(K_{1,n}^{\frac{1}{4}}) = 2n+2$, by Cases 1, 2 and 3.  \qed

\begin{theorem}\label{frac}
	For any $k\geq 3$, $\chi_d^t(G^{\frac{1}{k}}) \leq \chi_d^t(G^{\frac{1}{k+1}})$.
	\end{theorem}

\begin{figure}[h]
	\hspace{1.8cm}
	\begin{minipage}{4.5cm}
		\includegraphics[width=\textwidth]{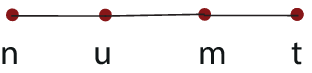}
	\end{minipage}
	\hspace{1cm}
	\begin{minipage}{6cm}
		\includegraphics[width=\textwidth]{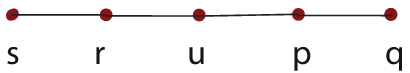}
	\end{minipage}
	\caption{\label{P4Conj} \small A part of a superedge in the proof of Theorem \ref{frac}, respectively.}
\end{figure}

\proof 
First we give a TD-coloring  to the vertices of $G^{\frac{1}{k+1}}$. Let  $P^{\{v,w\}}$ be an arbitrary superedge of $G^{\frac{1}{k+1}}$ with vertex set  $\{v, x^{\{v,w\}}_1, \ldots , x^{\{v,w\}}_{k}, w\}$. We have the following cases:
\begin{itemize}
	\item[Case 1)]
	There exists a vertex $u\in \{ x^{\{v,w\}}_1, \ldots , x^{\{v,w\}}_{k}\}$ such that other vertices of graph are not adjacent to all vertices with color class of vertex $u$. Consider the graph in Figure \ref{P4Conj}. Suppose that the vertex $u$ has the color $i$ and the vertex $n$ has color $\alpha$. The vertex $m$ is adjacent to  all vertices with color class $j$ and $j\neq i$ and  the vertex $n$ is adjacent to  all vertices with color class $k$ and $k\neq i$. Since $k\geq 3$, without loss of generality, suppose that $m\neq v$. We have two cases:

	\begin{itemize}
		\item[Case i)]
		The color of the vertex $m$ is not $\alpha$. We omit the vertex $u$   and put an edge between $n$ and $m$. So without adding a new color we have a TD-coloring  for this new graph.
		\item[Case ii)]
		The color of the vertex $m$ is $\alpha$. Since the vertex  $u$ is adjacent to  color class $\alpha$, so any other vertices does not  have color $\alpha$. In this case, by removing the vertex $m$ and putting   an edge between $u$ and $t$, we have a TD-coloring for this new graph. Because the vertex $t$ is adjacent to color class which is not $\alpha$, the color of $t$ is not $i$ (because of our assumptions), the vertex  $n$ is adjacent to  all vertices with color class $k$ and the vertex $u$ is adjacent to all vertices with color class $\alpha$.
	\end{itemize}
	
	\item[Case 2)]
	For every vertex $u\in \{ x^{\{v,w\}}_1, \ldots , x^{\{v,w\}}_{k}\}$, there exists a vertex such that is adjacent to all vertices with color of vertex $u$. Consider the graph in Figure \ref{P4Conj}. Suppose that the vertex $u$ has the color $i$ and the  vertex $p$ has color $j$ and the vertex $p$ is adjacent to all vertices with color $i$. We have two cases:
	\begin{itemize}
		\item[Case i)]
		The color of vertex $q$ is not $i$. We omit the vertex $r$   and put an edge between $u$ and $s$. So without adding a new color we have a TD-coloring  for this new graph, since there is no other vertex with color $i$.
		\item[Case ii)]
		The color of the vertex $q$ is  $i$. In this case the vertex $r$ is adjacent to color class of vertex $s$ and the color of vertex $s$ does not use for other vertices. Now we omit the vertex $u$ and put an edge between $r$ and $p$. If the color of $r$ is $j$, then we change it to $i$ and obviously since the vertex $s$ was adjacent to a color class except $j$, so we have a TD-coloring. If the color of the vertex $r$ is not $j$, then  we do not change its color and so we have a TD-coloring  again. 
	\end{itemize}
\end{itemize}
Now we give the same TD-coloring for  all superedges.  So we have the result.  \qed

 \begin{figure}
 	\begin{center}
 		\includegraphics[width=2in]{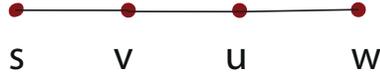}
 		\caption{A superedge in $G^{\frac{1}{3}}$.}
 		\label{P4}
 	\end{center}
 \end{figure}
 
 Here we prove that Theorem \ref{frac} is true for $k=2$, too.

 \begin{figure}
 	\begin{center}
 		\includegraphics[width=4in]{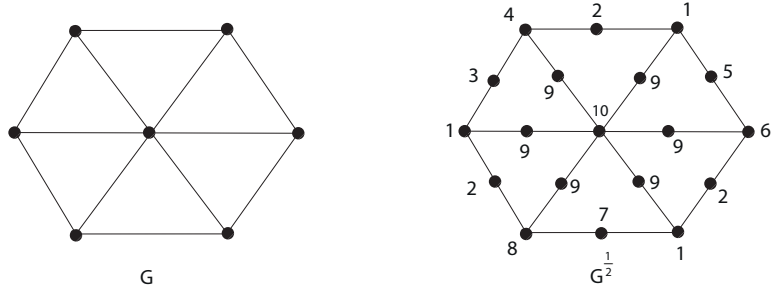}
 		\caption{An example which shows that Theorem \ref{last} is not true for $k=2$. }
 		\label{g12}
 	\end{center}
 \end{figure}

 \begin{theorem}
 For any graph $G$, 	$\chi_d^t(G^{\frac{1}{2}}) \leq \chi_d^t(G^{\frac{1}{3}}).$
  \end{theorem}
  \proof 
 First we give a TD-coloring  to the vertices of $G^{\frac{1}{3}}$. Let  $P^{\{s,w\}}$ be an arbitrary superedge of $G^{\frac{1}{3}}$ with vertex set  $\{s, x^{\{s,w\}}_1=v, x^{\{s,w\}}_{2}=u, w\}$ (see Figure \ref{P4}) and suppose that  the vertex $v$ has the color $\alpha$. We have the following cases:
 \begin{itemize}
 	\item[Case 1)]
 The vertices 	$u$ and $s$ are adjacent with a vertex with a color class which is not $\alpha$. we have two cases:
 	\begin{itemize}
 		\item[Case i)]
 		The color of vertices $u$ and $s$ are different.  In this case, by removing vertex $v$ and put an edge between $u$ and $s$, we have a TD-coloring for this new graph. Because  two vertices $u$ and $s$ are adjacent with a vertex with  color class which is not $\alpha$.
 		\item[Case ii)]
 		The color of vertices $u$ and $s$ are the same. Suppose that $u$ and $s$ have color $\beta$. In this case $\beta$ does not use for any other vertices. So $w$ is adjacent with a vertex with color class except $\beta$. Now we remove vertex $u$ and  put an edge between $v$ and $w$. So we have a TD-coloring for this new graph.
 	\end{itemize}
 	
 	\item[Case 2)]
 The vertex 	$u$ is adjacent to  all vertices with color class $\alpha$. We have two cases:
 	\begin{itemize}
 		\item[Case i)]
 		The color of the vertex $w$ is not $\alpha$. Suppose that the vertex $u$ has color $\gamma$.  If the vertex $v$ is adjacent with all vertices with color $\gamma$, and if the color of $s$ is  $\gamma$, we remove the vertex $u$ and  put an edge between $v$ and $w$. But  if the color of vertex $s$ is not $\gamma$, then we remove the vertex $u$, put an edge between $v$ and $w$ and give the color $\gamma$ to the vertex $w$. So we have a TD-coloring for this new graph. If the vertex $v$ is adjacent to  all vertices with color except $\gamma$ (vertex $s$), then we remove the vertex $u$  and  put an edge between $v$ and $w$. So we have  a TD-coloring for this new graph.
 		\item[Case ii)]
 		The color of the vertex  $w$ is $\alpha$. We have two new cases. First, the vertex $v$ is adjacent to a vertex with color class $\gamma$. Any adjacent vertex with $w$ is not adjacent to  vertex with color class $\alpha$ (except $u$). So we remove the vertex  $u$, put an edge between $v$ and $w$  and give the color $\gamma$ to $w$. This is a TD-coloring for this new graph.  Second, $v$ is not adjacent with color class $\gamma$. So the color of the vertex $s$ does not use any more. Also the vertex $s$ is not adjacent to vertex  with color class $\alpha$. So we remove vertex $v$ and put an edge between $s$ and $u$. This is a TD-coloring for this new graph.
 	\end{itemize}
 	
 	\item[Case 3)]
 The vertex	$s$ is adjacent to  all vertices with color class $\alpha$. We have two cases:
 	\begin{itemize}
 		\item[Case i)]
 		If $v$ is the only vertex which has color $\alpha$,  then we remove the vertex $u$ and put an edge between $v$ and $w$ when $v$ is adjacent with color class of vertex $s$, and remove  the vertex $v$ and put an edge between $s$ and $u$ and give the color $\alpha$ to $s$ when $v$ is adjacent with color class of vertex $u$. This is a TD-coloring for this new graph.
 		\item[Case ii)]
 		If there exist some vertices with color $\alpha$, then the vertex $u$ is adjacent with color class except $\alpha$. So we remove $v$ and put an edge between $s$ and $u$. This is a TD-coloring for this new graph.
 	\end{itemize}
 \end{itemize}
 We apply this TD-coloring for all superedges. So we obtain a TD-coloring for  $G^{\frac{1}{2}}$. Therefore we have $\chi_d^t(G^{\frac{1}{2}}) \leq \chi_d^t(G^{\frac{1}{3}})$.  \qed

 \begin{figure}
 	\begin{center}
 		\includegraphics[width=5.5in]{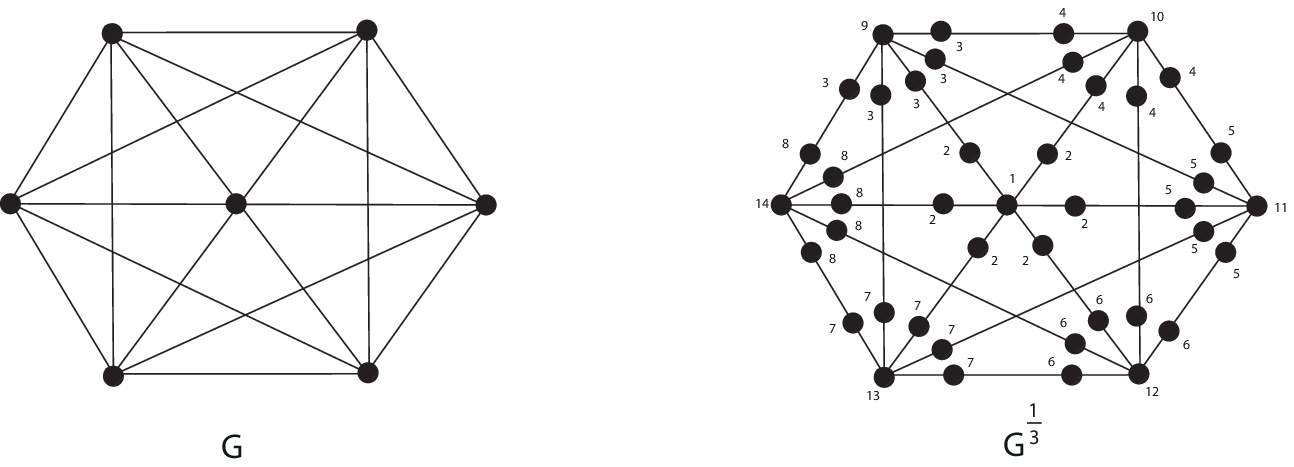}
 		\caption{\label{k=3} An example which show that Theorem \ref{last} is not true for $k=3$.}
 	\end{center}
 \end{figure}
 
 We end the paper with the following theorem:

 \begin{theorem}\label{last} 
 	If $G$ is  a graph with $m$ edges, then $\chi_d^t(G^{\frac{1}{k}}) \geq m$, for $k\geq 4$. 
 \end{theorem}
  \proof
 For $k=4$, in any superedge $P^{\{v,w\}}$ such as $\{v, x^{\{v,w\}}_1,  x^{\{v,w\}}_2, x^{\{v,w\}}_{3}, w\}$ without considering vertices $v$ and $w$,  there is a path graph of length two.  The vertex $x^{\{v,w\}}_2$ need to use a new color in at least one of its adjacent vertices, and we cannot use this color in any other superedges. 
 So we have the result. \qed

Note that Theorem \ref{last}  is not true for $k=2$. As an example,  for the graph $G$ of size $12$ in Figure \ref{g12} ,we have $\chi_d^t(G^{\frac{1}{2}})\ngeq 12$.
Also consider the graph $G$ of size $18$ in Figure \ref{k=3}. For this graph we have $\chi_d^t(G^{\frac{1}{3}})\ngeq 18$ which show that Theorem \ref{last}  is not true for $k=3$.

\end{document}